\documentstyle{article}
\def\R{\ifmmode{\rm I\mkern-3.1mu
R\mkern1mu}\else{\rm I\kern-.18em  R\hskip1pt\ 
}\fi\relax}  
\def\a{\alpha} 
\def\b{\beta} 

\def\G{\Gamma}

\def\t{\tau}
  
\def\th{\theta} 
\def\l{\lambda}

\def\s{\sigma}
\def\ch{\hat} 
\def\ep{\epsilon}

\def\so{\underline} 

\def\f{\rightarrow}
\def\q{\forall}
  
\def\v{\vdash}
 
\def\p{\succ}

\def\inc{\subseteq}  

\def\<{\ifmmode{\rm < \mkern-9.1mu
<\mkern0.3mu}\else{\rm I\kern-.16em
N\hskip0.5pt\ }\fi\relax} 

\def\>{\ifmmode{\rm > \mkern-9.1mu
>\mkern0.3mu}\else{\rm I\kern-.16em
N\hskip0.5pt\ }\fi\relax}

\begin{document}

\begin{center}
{\LARGE\bf  Une r\'eponse n\'egative 
\`a la conjecture de E.Tronci
pour les syst\'emes num\'eriques typ\'es}\\ [0,5cm]
\end{center}

\begin{center} 
\bf 
Karim NOUR\\
\rm
LAMA - Equipe de Logique\\
Universit\'e de Savoie\\
73376 Le Bourget du Lac\\
e-mail nour@univ-savoie.fr\\[1cm]
\end{center}

{\bf R\'esum\'e} {\it Un syst\`eme num\'erique est une suite de $\l$-termes normaux clos
distincts pour laquelle il existe des $\l$-termes clos pour les fonctions successeur et test
\`a z\'ero. Un syst\`eme num\'erique est dit ad\'equat ssi il existe un $\l$-terme clos 
pour la fonction  pr\'ed\'ecesseur. Un op\'erateur de mise en m\'emoire pour un  syst\`eme
num\'erique est un $\l$-terme clos qui simule ``l'appel-par-valeur'' dans le cadre de
``l'appel-par-nom''. E. Tronci a conjectur\'e le r\'esultat suivant : un syst\`eme num\'erique
est ad\'equat s'il poss\`ede un op\'erateur de mise en m\'emoire. Nous donnons, dans cet
article, une r\'eponse n\'egative \`a la conjecture de E. Tronci mais uniquement pour les
syst\`emes num\'eriques typable dans le syst\`eme $\cal F$. La conjecture de E. Tronci reste
sans solution en $\l$-calcul pur.}\\

{\bf Mots cl\'es :} Syst\`eme num\'erique ; $\l$-calcul ; Successeur ; Test \`a z\'ero ;
Pr\'ed\'ecesseur ; Op\'erateur de mise en m\'emoire ; Syst\`eme num\'erique ad\'equat ;
Appel-par-valeur ; Appel-par-nom ; Syst\`eme $\cal F$.\\

{\bf Abstract} {\it A numeral system is a sequence of an infinite different closed normal
$\l$-terms which has closed $\l$-terms for successor and  zero test. A numeral system is
said  adequate iff it has a closed $\l$-term for predecessor. A storage operator for a 
numeral system is a closed $\l$-term which simulate ``call-by-value'' in the context of a
``call-by-name'' strategy. E. Tronci conjectured the following result :  a numeral system is
adequate if it has a storage operator. This paper gives a negative answer to this
conjecture for the numeral systems typable in the J.-Y. Girard type system $\cal F$. The
E. Tronci's conjecture remains open in pure $\l$-calculus.} \\

{\bf Keywords :} Numeral system ; $\l$-calculus ; Successor ; Zero test ; Predecessor ;
Storage operator ; Adequate numeral system ; Call-by-value ; Call-by-name ;
Type system $\cal F$.\\

\section{Introduction}
Un syst\`eme num\'erique est une suite de $\l$-termes normaux clos distincts $\bf
d$ \rm = $d_0 , d_1 ,..., d_n , ...$ pour laquelle il existe des $\l$-termes clos $S_d$ et 
$Z_d$ pour les fonctions successeur et test \`a z\'ero. Un syst\`eme
num\'erique est dit ad\'equat ssi il existe un $\l$-terme clos $P_d$ pour la fonction 
pr\'ed\'ecesseur. H. Barendregt a d\'emontr\'e dans [1] qu'un syst\`eme num\'erique est
ad\'equat ssi toutes les fonctions r\'ecursives totales sont repr\'esentables dans le
syst\`eme. \\

La diff\'erence entre notre d\'efinition d'un syst\`eme num\'erique et celle propos\'ee par
H. Barendregt (voir [1]) est le fait d'imposer aux $\l$-termes $d_i$ d'\^etre normaux
et distinctes. En effet ces conditions permettent, pour des strat\'egies de r\'eduction
gagnantes, de trouver la valeur exacte d'une fonction num\'erique totale calcul\'ee sur des
entiers.\\

Une des strat\'egies de r\'eduction gagnantes est la r\'eduction gauche (it\'eration de la
r\'eduction de t\^ete not\'ee $\p$). Mais pour cette strat\'egie l'argument d'une fonction est
calul\'e le nombre de fois o\`u la fonction l'utilise. Les op\'erateurs de mise en m\'emoire
ont \'et\'e introduits par J.-L. Krivine pour rem\'edier \`a ce d\'efaut. \\

Un $\l$-terme clos $O_d$ est dit op\'erateur de mise en m\'emoire pour un syst\`eme
num\'erique $\bf d$ ssi pour tout $n \in {\bf N}$, il existe un $\l$-terme clos $\t_n
\simeq_{\b} d_n$ tel que pour tout $\th_n \simeq\sb{\b} d_n$,  $(O_d ~ \th_n ~ f) \p (f ~
\t_n)$ (o\`u $f$ est une nouvelle variable).\\

Nous allons justifier cette d\'efinition. Soit $F$ un $\l$-terme (pour une fonction), et
$\th_n$ un $\l$-terme  $\b$-\'equivalent \`a $d_n$. Durant la r\'eduction gauche de $(F ~
\th_n)$, $\th_n$ sera r\'eduit chaque fois qu'il arrive en t\^ete. Au lieu de r\'eduire $(F ~
\th_n)$, effectuons la r\'eduction de t\^ete de $(O_d ~ \th_n ~ F)$. La r\'eduction de $ (O_d
~ \th_n ~ F) = \{(O_d ~ \th_n ~ f)\}[F/f]$ commence par amener $(O_d ~ \th_n ~ f)$ \`a sa
forme normale de t\^ete qui est $(f ~ \t_n)$, et puis r\'eduire $(F ~ \t_n)$. Dans la
r\'eduction de $(O_d ~ \th_n ~ F)$, $\th_n$ est calcul\'e le premier, et le r\'esultat est
donn\'e \`a $F$ comme argument. $O_d$ a donc mis en m\'emoire le r\'esultat $\t_n$, avant de
le donner \`a la fonction $F$. Donc la r\'eduction de t\^ete $(O_d ~ \th_n ~ F) \p
(F ~ \t_n)$ d\'epend seulement de $\th_n$ et pas de $F$.\\ 

J.-L. Krivine a d\'emontr\'e dans [4] que, dans le syst\`eme de typage $\cal F$ de
J.-Y. Girard, le type  $N$*$ \f \neg\neg N$  convient pour les op\'erateurs de mise en
m\'emoire pour le syst\`eme num\'erique de Church : o\`u $N$ est le type des entiers de
Church, et l'op\'eration $*$ est la simple traduction de G\H{o}del qui associe \`a chaque
formule $F$ la formule $F$* obtenue en remplacant dans $F$ chaque variable de type par sa
n\'egation.\\

Nous d\'emontrons dans ce papier que chaque syst\`eme num\'erique ad\'equat poss\`ede un
op\'erateur de mise en m\'emoire. E. Tronci a conjectur\'e qu'un syst\`eme num\'erique est
ad\'equat s'il poss\`ede un op\'erateur de mise en m\'emoire.\\

Nous donnons, ensuite, une r\'eponse n\'egative \`a la conjecture de E. Tronci mais uniquement
pour les syst\`emes num\'eriques typable dans le syst\`eme $\cal F$. Nous construisons
donc un type clos $E$, une suite de $\l$-termes normaux clos distincts ${\bf e} = e_0 , e_1
,..., e_n , ...$, et des  $\l$-termes clos $S_e$, $Z_e$, et $O_e$ tels que :

\begin{itemize} 
\item[] -- Si $t$ est un $\l$-terme normal clos, alors $\v_{\cal F} t : E$ ssi $t = e_i$ o\`u
$i \in {\bf N}$.
\item[] -- $\v_{\cal F} S_e : E \f E$ et $(S_e ~ e_n) \simeq_{\b} e_{n+1}$ pour tout $n \in
{\bf N}$. 
\item[] -- $\v_{\cal F} Z_e : E \f B$ ($B$ est le type des
Bool\'eens du syst\`eme $\cal F$), $(Z_e ~ e_0) \simeq\sb{\b} \l x \l y x$ et
$(Z_e  ~ e_{n+1}) \simeq\sb{\b} \l x \l y  y$ pour tout $n \in {\bf N}$.
\item[] -- $\v_{\cal F} O_e : E$*$ \f \neg\neg E$, et, pour tout $n \in {\bf N}$, il existe un
$\l$-terme clos $\t_n \simeq_{\b} e_n$ tel que pour tout $\th_n \simeq\sb{\b} e_n$, 
$(O_d ~ \th_n ~ f) \p (f ~ \t_n)$. 
\item[] -- Il n'existe pas un $\l$-terme clos $P_e$ tel que $\v_{\cal F} P_e : E \f E$ et
$(P_e ~ e_{n+1}) \simeq_{\b} e_n$ pour tout $n \in {\bf N}$. 
\end{itemize}

La conjecture de E. Tronci reste sans solution en $\l$-calcul pur.

\section{Notations et d\'efinitions}

\subsection{Le $\l$-calcul pur}

{\bf Notations} :
\begin{itemize}
\item[] 1) La $\b$-\'equivalence est not\'ee $u \simeq\sb{\b} v$. 
\item[] 2) Si $u$ et $v$ sont deux $\l$-termes, alors on note $<u,v>$ le $\l$-terme $\l x (x ~
u~v)$.
\item[] 3) On note $T$ (pour True) le $\l$-terme $\l x \l y x$ et $F$ (pour False) le
$\l$-terme $\l x \l y y$. 
\item[] 4) Pour tous $\l$-termes $u,v$, on d\'efinit $(u^n ~ v)$ par induction : $(u^0 ~ v) =
v$ et $(u^{n+1} ~ v) = (u ~ (u^n ~ v))$. Pour chaque entier $n$, on d\'efinit {\it l'entier de
Church} $\so{n} = \l x\l f(f^n ~ x)$. 
\item[] 5) La notation $\s(t)$ repr\'esent le r\'esultat d'une substitution simultan\'ee $\s$
sur les variables libres de $t$ apr\`es un r\'enommage de ses variables li\'ees. 
\item[] 6) On note $\Theta = (U ~ U)$ o\`u $U = \l x \l f (f ~ (x ~ x ~ f))$. Le $\l$-terme
$\Theta$ est appel\'e {\it le point fixe de Turing}. 
\end{itemize}

{\bf D\'efinitions} :  Un $\l$-terme $t$ soit il poss\`ede un {\it redex de t\^ete} [i.e. $t=\l
x_1 ...\l x_n (\l x u ~ v ~ v_1 ... v_m)$, le redex de t\^ete est $(\l x u ~ v)$], soit il est
en {\it forme normale de t\^ete} [i.e. $t=\l x_1 ...\l x_n (x ~ v_1 ... v_m)$]. La notation
$u \p v$ signifie que $v$ est obtenue \`a partir de $u$ apr\`es quelques pas de r\'eductions
de t\^ete.  Un $\l$-terme est dit {\it r\'esoluble} si sa r\'eduction de t\^ete
termine.  \\

Les r\'esultats suivants sont bien connus (voir [3] et [4]). \\

{\bf Th\'eor\`eme  1} \\
{\it 1) Si $t$ est $\b$-\'equivalent \`a une forme normale de t\^ete, alors $t$
est r\'esoluble.} \\
{\it 2) Si $u \p v$, alors, pour toute substitution $\s$, $\s(u) \p \s(v)$.}\\
{\it 3) Si $u \p v$, alors, pour toute suite $w_1,...,w_n$, il existe un
$\l$-terme $w$ tel que $(u ~ w_1 ... w_n) \p w$ et $(v ~ w_1 ... w_n) \p w$.} \\

{\bf D\'efinition} : On d\'efinit sur les $\l$-termes une relation d'\'equivalence  $\sim$ par
: $u \sim v$ ssi il existe un $\l$-terme $t$, tel que $u \p t$, et $v \p t$. \\

Donc, si $t$ est r\'esoluble, alors $u \sim t$ ssi $u$ est r\'esoluble, et poss\`ede
la m\^eme forme normale de t\^ete que $t$. Si $u$ est une forme normale de t\^ete, alors $t
\sim u$ signifie que $u$ est la forme normale de t\^ete de $t$. \\

D'apr\`es le th\'eor\`eme 1, on obtient les r\'esultats suivants (voir [4]). \\

{\bf Th\'eor\`eme  2} \\
{\it 1) Si $u \sim v$, alors, pour toute substitution $\s$, $\s(u) \sim \s(v)$.}\\
{\it 2) Si $u \sim v$, alors, pour toute suite $w_1,...,w_n$, $(u ~ w_1 ... w_n)  \sim (v ~
w_1 ... w_n)$.} 

\subsection{Le syst\`eme {$\cal F$}}

{\bf D\'efinition} : Les types du syst\`eme {$\cal F$} sont construits \`a partir des variables
de type $X,Y,Z,...$ et une constante $\perp$ (pour l'absurde) en utilisant les op\'erations
suivantes :  
\begin{itemize}
\item[] -- Si $U$ et $V$ sont des types, alors $U \f V$ est un type.
\item[] -- Si $V$ est un type, et $X$ est une variable de type, alors $\q X V$ est un type.\\
\end{itemize}

On d\'efinit d'une mani\`ere usuelle les {\it variables libres} et les {\it variables li\'ees}
d'un type. \\

{\bf D\'efinition} : Soient $t$ un $\l$-terme, $A$ un type, et $\G = x_1 : A_1 ,..., x_n : A_n$
un contexte. On d\'efinit par les r\`egles suivantes la notion ``$t$ est de type $A$ dans
$\G$'' ; cette notion est not\'ee  $\G\v_{\cal F} t:A$.   

\begin{center} 
(1) $\G\v_{\cal F} x_i:A_i$ $(1\leq i\leq n)$
\end{center} 
\begin{minipage}[t]{210pt}
$  ( 2 ) \quad \displaystyle\frac{ \G,x:A \v_{\cal F} t:B } { \G\v_{\cal F} \l xt:A \f B }$ \\
\end{minipage} 
\begin{minipage}[t]{210pt}\sl 
$  ( 3 ) \quad  \displaystyle\frac{ \G\v_{\cal F} u:A \f B \quad \G\v_{\cal F} v:A} { \G\v_{\cal F} (u)v:B
}$ \\ \end{minipage}
\begin{minipage}[t]{210pt}
$  ( 4 ) \quad \displaystyle\frac{ \G\v_{\cal F} t:A } { \G\v_{\cal F} t:\q XA }$  {\rm (*)}\\
\end{minipage} 
\begin{minipage}[t]{210pt}\sl 
$  ( 5 ) \quad  \displaystyle\frac{ \G\v_{\cal F} t:\q XA} { \G\v_{\cal F} t:A[G/X] }$  {\rm (**)}\\ 
\end{minipage}

Avec les conditions suivantes : \\
(*) $X$ n'est pas libre dans $\G$.\\
(**) $G$ est un type. \\

Le syst\`eme $\cal F$ poss\`ede les propri\'et\'es suivantes (voir [2] et [3]).\\

{\bf Th\'eor\`eme  3} \\{\it 
1) Un type est pr\'eserv\'e durant une $\b$-r\'eduction.\\
2) Un $\l$-terme typable est fortement normalisable.}\\

Le type $F_1 \f (F_2 \f(...\f (F_n \f G)...))$ est not\'e $F_1,F_2,...,F_n \f
G$ et le type $F \f \perp$ est not\'e $\neg F$.\\

Les lemmes 1 et 2 seront tr\`es utiles pour nos d\'emonstrations. Le lemme 1 (resp. le lemme
2) a \'et\'e d\'emontr\'e dans [5] (resp. dans [2] et [3]). \\

{\bf Lemme 1} \\
{\it 1) Soit $X$ une variable de type. Si $\G\v_{\cal F} t:X$, alors $t$ ne commence pas par
$\l$.\\ 2) Si $\G\v_{\cal F} \l x t:A \f B$, alors $\G , x : A \v_{\cal F} t: B$. \\
3) Si $\G , x : A \f B \v_{\cal F} (x ~ u_1 ... u_n) :C$, alors $\G , x : A \f B \v_{\cal F} 
u_1 : A$. \\
4)  Si $\G , x : A_1,...,A_m \f X \v_{\cal F} (x ~ u_1 ... u_n) :C$, alors $X$ est libre dans
$C$.} \\

{\bf Lemme 2} {\it Si $x_1 : A_1,...,x_n : A_n \v_{\cal F} t : A$, alors, pour toute
variable $X$ et tout type $G$, $x_1 : A_1[G/X],...,x_n : A_n[G/X] \v_{\cal F} t : A[G/X]$.}\\

Dans le syst\`eme  ${\cal F}$ on a la possibilit\'e de d\'efinir les types de donn\'ees. Soit
$B = \q X \{ X,X \f X \}$ ({\it le type des Bool\'eens}) et $N = \q X \{ X,(X \f X) \f X \}$
({\it le type des entiers}). On a les r\'esultats suivants (voir [2] et [3]).\\

{\bf Th\'eor\`eme  4} {\it Soit $t$ un $\l$-terme normal clos. \\ 
1) $\v_{\cal F} t : B$ ssi $t = T$ ou $t = F$.\\ 
2) $\v_{\cal F} t : N$ ssi il existe $n \in {\bf N}$ tel que $t = \so{n}$.}\\ 

{\bf D\'efinition} : On note encore $\cal F$ le syst\`eme logique sousjacent au syst\`eme de
typage $\cal F$, et on \'ecrit $\G \v_{\cal F} A$ si $A$ est d\'emontrable \`a partir des
formules de $\G$ en utilisant les r\`egles du syst\`eme logique ${\cal F}$. \\

Il est claire que : $A_1,...,A_n \v_{\cal F} A$ ssi il existe un $\l$-terme $t$ tel que 
$x_1 : A_1,...,x_n : A_n \v_{\cal F} t : A$. Ce r\'esultat est connu sous le nom de ``la
corresponce du Curry-Howard''. 

\subsection{Le syst\`eme ${\cal F_C}$ et la traduction de G\H{o}del}

{\bf D\'efinition} : On ajoute au syst\`eme logique $\cal F$ la r\`egle : 
\begin{center}
$  ( 0 ) \quad \displaystyle\frac{ \G \v \neg \neg A } { \G \v A }$
\end{center}
Cette r\`egle axiomatise la logique classique au dessus de la logique intuitionniste.\\
On note $\cal F_C$ ce nouveau syst\`eme et on \'ecrit $\G \v_{\cal F_C} A$ si $A$ est
d\'emontrable \`a partir des formules de $\G$ dans le syst\`eme $\cal F_C$. On a le
r\'esultat suivant (voir [2]).\\

{\bf Th\'eor\`eme  5} {\it Le syst\`eme $\cal F_C$ est non contradictoire (i.e. $\not \v_{\cal
F_C} \q XX$).}\\

{\bf D\'efinition} : Pour chaque formule $A$ de $\cal F_C$, on d\'efinit la formule $A$* par : 
\begin{itemize}  
\item[] - Si $A=\perp$, alors  $A$*$= A$ ;
\item[] - Si $A=X$, alors  $A$*$ = \neg X$ ;
\item[] - Si $A=B \f C$, alors $A$*$ = B$*$ \f C$* ;
\item[] - Si $A=\q X B$, alors $A$*$ = \q X B$*.
\end{itemize}
$A$* est appel\'ee {\it la traduction de G\H{o}del} de $A$. \\

On a le r\'esultat suivant (voir [3]).\\

{\bf Th\'eor\`eme  6} {\it Si $\v_{\cal F_C} A$, alors $\v_{\cal F} A$*.}

\section{Les syst\`emes num\'eriques}

\subsection{Les syst\`emes num\'eriques en $\l$-calcul pur}Ê

{\bf D\'efinition} : Un {\it syst\`eme num\'erique} est une suite de $\l$-termes normaux clos
distincts $\bf d$ \rm = $d_0 , d_1 ,..., d_n , ...$ pour laquelle il existe des $\l$-termes
clos $S_d$ et $Z_d$ tels que : 
\begin{center} $(S_d ~ d_n) \simeq\sb{\b} d_{n+1}$ pour tout $n \in {\bf N}$ \\
et\\ 
$(Z_d ~ d_0) \simeq\sb{\b} T$ \\
$(Z_d  ~ d_{n+1}) \simeq\sb{\b} F$ pour tout $n \in {\bf N}$ 
\end{center}
Les $\l$-termes $S_d$ et $Z_d$ sont appel\'es {\it successeur} et {\it test \`a z\'ero} pour 
$\bf d$.\\

Chaque syst\`eme num\'erique peut \^etre consid\'erer comme un codage des entiers en 
$\l$-calcul et donc on peut repr\'esenter les fonctions num\'eriques totales de la mani\`ere
suivante.\\

{\bf D\'efinition} : Une fonction num\'erique totale $\phi : {\bf N}^p \f {\bf N}$ is dite {\it
$\l$-d\'efinissable dans le  syst\`eme num\'erique {\bf d}} ssi il existe un $\l$-terme
$F_{\phi}$ tel que pour tout $n_1,...,n_p \in {\bf N}$
\begin{center}
 $(F_{\phi} ~ d_{n_1}... d_{n_p}) \simeq\sb{\b} d_{\phi(n_1,...,n_p)}$ 
\end{center}

{\bf D\'efinition} : Un syst\`eme num\'erique {\bf d} is dit {\it ad\'equat} ssi il existe un
$\l$-terme clos $P_d$ tel que 
\begin{center} 
$(P_d ~ d_{n+1}) \simeq\sb{\b} d_n$ pour tout $n \in {\bf N}$. 
\end{center}
Le $\l$-terme $P_d$ est appel\'e {\it pr\'ed\'ecesseur} pour {\bf d}. \\

H. Barendregt a d\'emontr\'e que (voir [1]) :\\

{\bf Th\'eor\`eme 7} {\it Un syst\`eme num\'erique {\bf d} \it est ad\'equat ssi toutes les
fonctions num\'eriques r\'ecursives totales sont $\l$-d\'efinissables dans} {\bf d}. \\ 

{\bf Exemples} : \\
1) Un exemple simple d'un syst\`eme num\'erique ad\'equat est le {\it
syst\`eme num\'erique de Church} {\bf \so{n}} = $\so{0},\so{1},...,\so{n},...$. Il est facile
de v\'erifier que : \begin{itemize}
\item[] $\so{S}=\l n \l x \l f(f ~ (n ~ f ~ x))$,
\item[] $\so{Z}=\l n (n ~  T ~ \l x F)$,   
\item[] $\so{P}=\l n(n ~ U ~ <\so{0},\so{0}> ~ T)$ o\`u $U=\l x<(\so{S}~ (x ~ T)),(x ~ F)>$. 
 \end{itemize}
sont des $\l$-termes pour le successeur, le test \`a z\'ero, et le pr\'ed\'ecesseur pour {\bf
\so{n}}. \\ 
2) Nous avons donn\'e dans [6] un exemple d'un syst\`eme num\'erique non ad\'equat. \hfill
$\Box$ \\

{\bf D\'efinition} : Soient {\bf d} un syst\`eme num\'erique et $O_d$ un $\l$-terme clos. On
dit que $O_d$ est un {\it op\'erateur de mise en m\'emoire} pour {\bf d} ssi pour tout $n \in
{\bf N}$, il existe un $\l$-terme clos $\t_n \simeq\sb{\b} d_n$, tel que, pour tout $\th_n
\simeq\sb{\b} d_n$,  $(O_d ~ \th_n f) \p (f ~ \t_n)$ o\`u $f$ est une nouvelle variable.  \\

{\bf Exemple} : Soit $O_N = \l n\l f(n ~ f ~ J ~ \so{0})$
o\`u $J = \l x\l y(x ~ (\so{S} ~ y))$. Il est facile de v\'erifier que pour tout $\th_n
\simeq\sb{\b} \so{n}$, $(O_N ~ \th_n ~ f) \p (f ~ (\so{S}^n ~ \so{0}))$. Donc $O_N$ est un
op\'erateur de mise en m\'emoire pour {\bf \so{n}}. \hfill $\Box$ \\

{\bf Remarque} : J.-L. Krivine autorise, dans sa d\'efinition des op\'erateurs de mise en
m\'emoire, le $\l$-terme $\t_n$ de contenir des variables libres qui peuvent \^etre
remplac\'ees par des $\l$-termes qui ne d\'ependent que de $\th_n$. Avec cette d\'efinition on
garde aussi tous les r\'esultats de ce papier. \hfill $\Box$ \\

{\bf Th\'eor\`eme  8} {\it Chaque syst\`eme num\'erique ad\'equat poss\`ede un op\'erateur de
mise en m\'emoire.} \\

{\bf Preuve} Soit {\bf d} un syst\`eme num\'erique ad\'equat. \\ 
Soit $O_d = (\Theta ~ H_d)$ o\`u
$H_d = \l h \l n \l f ((Z_d ~ n) ~ (f ~ d_0) ~ (h ~ (P_d ~ n) ~ \l x (f~ (S_d ~ x))))$.\\ 
D\'emontrons (par r\'ecurrence sur $i$) que, pour tout $i \in {\bf N}$ et pour tout
$\th_i \simeq\sb{\b} d_i$, $(O_d ~ \th_i ~ f) \sim (f ~ ({S_d}^i ~ d_0))$.  
\begin{itemize}
\item Pour $i=0$, 
\begin{eqnarray} 
(O_d ~ \th_0 ~ f) &\sim &(H_d ~ (\Theta ~ H_d) ~ \th_0 ~ f)  \nonumber \\
&\sim &((Z_d ~ \th_0) ~ (f ~ d_0) ((\Theta ~ H_d) ~ (P_d ~ \th_0) ~ \l x (f ~ (S_d ~ x)))) 
\nonumber 
\end{eqnarray}
 Comme $(Z_d ~ d_0) \simeq\sb{\b} T$, alors
$(Z_d ~ \th_0) \p T$, et, d'apr\`es le th\'eor\`eme 1, $(O_d ~ \th_0 ~ f) \sim (f ~ d_0)$. 
\item Supposons le r\'esultat vrai pour $i$, et prouvons le pour $i+1$. 
\begin{eqnarray} 
(O_d ~ \th_{i+1} ~ f) &\sim &(H_d ~  (\Theta ~ H_d) ~ \th_{i+1} ~ f)  \nonumber \\
&\sim &((Z_d ~ \th_{i+1}) ~ (f ~ d_0) ~ ((\Theta ~ H_d) ~ (P_d ~ \th_{i+1})) ~ \l x (f ~ (S_d ~
x)))) \nonumber 
\end{eqnarray}
Comme $(Z_d ~ d_{i+1}) \simeq \sb{\b} F$, alors $(Z_d ~ \th_{i+1}) \p F$, et,
d'apr\`es le th\'eor\`eme  1, \\ $(O_d ~ \th_{i+1} ~ f) \sim (O_d ~ (P_d ~ \th_{i+1}) ~ \l x
(f ~ (S_d ~ x)))$. Mais $(P_d ~ \th_{i+1}) \simeq\sb{\b} d_{i}$, alors, par hypoth\`ese
d'induction, $(O_d ~ (P_d ~ \th_{i+1}) ~ f) \sim (f ~ ({S_d}^i d_0))$, et  
\begin{eqnarray}  
(O_d ~ (P_d ~ \th_{i+1}) ~ \l x (f ~ (S_d ~ x))) &\sim &(\l x (f ~ (S_d ~ x)) ~ ({S_d}^i d_0))
\nonumber \\
&\sim &(f ~ ({S_d}^{i+1} d_0) )\nonumber 
\end{eqnarray}
 D'o\`u, pour tout $i \in {\bf N}$ et pour tout $\th_i \simeq_{\b} d_i$, $(O_d ~ \th_i ~ f) \p
(f ~ ({S_d}^{i} d_0))$.  \hfill $\Box$ \\ \end{itemize}

E. Tronci a conjectur\'e le r\'esultat suivant :\\

{\bf Conjecture} {\it Un syst\`eme num\'erique est ad\'equat s'il poss\`ede un op\'erateur de
mise en m\'emoire.}\\

Nous donnons dans ce papier une r\'eponse n\'egative \`a cette conjecture mais uniquement
pour les syst\`emes num\'eriques typable dans le syst\`eme ${\cal F}$.

\subsection{Les syst\`emes num\'eriques typ\'es}

{\bf D\'efinition} : Un {\it syst\`eme num\'erique typ\'e} est une paire ${\cal D} = <D,{\bf
d}>$ o\`u $D$ est un type clos du syst\`eme ${\cal F}$, et $\bf d$ \rm = $d_0 , d_1
,..., d_n , ...$ est une suite de $\l$-termes normaux clos tels que :  
\begin{itemize} \item[] --  Si $t$ est un $\l$-terme normal, alors $\v_{\cal F} t:D$ ssi
il existe $i \in {\bf N}$ tel que $t = d_i$. 
\item[] -- Il existe des $\l$-termes clos $S_d$ et $Z_d$ tels que:  
\begin{itemize}  
\item[]  * $\v_{\cal F} S_d : D \f D$ et $(S_d ~ d_n) \simeq\sb{\b} d_{n+1}$ pour tout $n \in
{\bf N}$ ;  
\item[]  * $\v_{\cal F} Z_d : D \f B$ et 
$(Z_d ~ d_n) \simeq\sb{\b} \cases { T &si $n=0$ \cr F &si $n \geq 1$ \cr}$.  
\end{itemize}
\end{itemize}
Les $\l$-termes $S_d$ et $Z_d$ sont appel\'es {\it successeur} et {\it test \`a z\'ero} pour 
$\cal D$.\\
 
{\bf D\'efinition} :  Un syst\`eme num\'erique typ\'e $\cal D$ est dite {\it ad\'equat} ssi il existe
un $\l$-terme clos $P_d$ tel que $\v_{\cal F} P_d : D \f D$ et $(P_d ~ d_{n+1}) \simeq\sb{\b}
d_n$ pour tout $n\in {\bf N}$. Le $\l$-terme $P_d$ est appel\'e {\it pr\'ed\'ecesseur} pour
$\cal D$.\\

{\bf Exemple} : Il est facile de v\'erifier que ${\cal N} = <N,{\bf \so{n}}>$ est un syst\`eme 
num\'erique typ\'e. \hfill $\Box$\\

{\bf D\'efinitions} :\\
1) Soient $D,E$ deux types clos. On dit que $D \inc E$ ssi pour tout
$\l$-terme clos $t$, si $\v_{\cal F} t : D$, alors $\v_{\cal F} t : E$. \\
2) Soit ${\cal D} = <D,{\bf d}>$ un syst\`eme num\'erique typ\'e tel que $D \inc D$*. Soit
$O_d$ un $\l$-terme clos. On dit que $O_d$ est un {\it op\'erateur de mise en m\'emoire} pour
$\cal D$ ssi $\v_{\cal F} O_d : D$*$ \f \neg \neg D$, et pour tout $n \in {\bf N}$, il existe
un $\l$-terme clos $\t_n \simeq\sb{\b} d_n$ et $\v_{\cal F} \t_n : D$ tel que, pour tout $\th_n
\simeq\sb{\b} d_n$,  $(O_d ~ \th_n f) \p (f ~ \t_n)$ o\`u $f$ est une nouvelle variable. \\

{\bf Exemple} : On peut v\'erifier que $N \inc N$* et $\v_{\cal F} O_N : N$*$\f\neg\neg
N$. Donc $O_N$ est un op\'erateur de mise en m\'emoire pour ${\cal N}$. \hfill $\Box$

\section{Le contre exemple}

Soit $P = \q X \q Y \{ ((X \f Y) \f X) \f X \}$ ($P$ est la loi de Pierce) et $Q = P \f \q
XX$. \\

{\bf Lemme 3} \\
1) {\it $\v_{\cal F_C} P$. \\
2) Il existe un $\l$-terme clos $t_P$ tel que $\v_{\cal F} t_P : P$*.} \\

{\bf Preuve}  \\
1) C'est un r\'esultat connu. Faisons la d\'emonstration.
\begin{eqnarray}
\neg X,  X, \neg Y \v_{\cal F_C}  \perp  &\Longrightarrow 
&\neg X,  X \v_{\cal F_C}  \neg \neg Y \nonumber \\
&\Longrightarrow &\neg X,  X \v_{\cal F_C} Y \nonumber \\
&\Longrightarrow &\neg X,  (X \f Y) \f X \v_{\cal F_C} X \f Y  \nonumber \\
&\Longrightarrow &\neg X,  (X \f Y) \f X \v_{\cal F_C} X  \nonumber \\
&\Longrightarrow &\neg X,  (X \f Y) \f X \v_{\cal F_C} \perp  \nonumber \\
&\Longrightarrow &(X \f Y) \f X \v_{\cal F_C} \neg \neg X  \nonumber \\
&\Longrightarrow &(X \f Y) \f X \v_{\cal F_C} X  \nonumber \\
&\Longrightarrow &\v_{\cal F_C} P.  \nonumber
\end{eqnarray}
2) D'apr\`es 1) et le th\'eor\`eme 6, on a $\v_{\cal F} P$*, donc il existe un
$\l$-terme $t_P$ tel que  $\v_{\cal F} t_P : P$*. Un exemple d'un tel $\l$-terme est $t_P =
\l x \l y (x ~ \l z \l \a (z ~ y) ~ y)$. En effet :
\begin{eqnarray}
z : \neg X, y :  X, \a :  Y \v_{\cal F} (z ~ y) : \perp  &\Longrightarrow 
&y :  X \v_{\cal F} \l z \l \a (z ~ y) :  \neg X \f \neg Y \nonumber \\
&\Longrightarrow &x : ( \neg X \f  \neg Y) \f  \neg X ,y :  X  \v_{\cal F} (x ~
\l z \l \a (z ~ y) ~ y) : \perp \nonumber \\
&\Longrightarrow &\v_{\cal F} t_p : P^{*}.  \nonumber 
\end{eqnarray}
\hfill $\Box$

{\bf Lemme 4} \\
1) {\it $\not \v_{\cal F} Q \f P$.\\
2) $\not \v_{\cal F_C} (Q \f P) \f Q$.} \\

{\bf Preuve}  \\
1) Un contexte $\G$ est dit {\it bon} ssi $\G$ est de la forme $[ \a : Q , \{ x_i : (X_i \f
Y_i) \f X_i \}_{1\leq i \leq n} , \{ y_j : X_j\}_{1\leq j \leq m} ]$ o\`u :
\begin{itemize}
\item[] -- $X_i \neq X_j$  ($1\leq i < j \leq n$) ;
\item[] --  $Y_i \neq Y_j$ ($1\leq i < j \leq m$) ; 
\item[] --  $X_i \neq Y_j$ ($1\leq i \leq n$) et ($1\leq j \leq m$).
\end{itemize}
Il suffit de d\'emontrer que pour tout contexte bon $\G$ il n'existe pas de $\l$-terme $u$
tel que $\G \v_{\cal F} u : P$. Nous d\'emontrons ceci par induction sur $u$. \\
$u$ ne peut pas \^etre une variable. Si $u = (z ~ u_1...u_m)$ ($m \geq 1$), alors $z = \a$,
et $\G \v_{\cal F} u_1 : P$. Ce qui est impossible par hypoth\`ese d'induction. Donc $u =
\l x v$, et $\G , x : (X \f Y) \f X \v_{\cal F} v : X$ o\`u $X,Y$ sont des nouvelles
variables diff\'erentes. $v$ ne peut pas \^etre ni une variable ni un $\l$-terme qui commence
par $\l$. Donc $v = (z ~ v_1...v_m)$ ($m \geq 1$) et $z \not = x_i, y_j$. Il reste,
donc, deux cas \`a voir :
\begin{itemize}
\item[] -- Si $z = \a$, alors $\G , x : (X \f Y) \f X \v_{\cal F} v_1 : P$. Ce qui est
impossible par hypoth\`ese d'induction. 
\item[] -- Si $z = x$, alors $m = 1$, et $\G , x : (X
\f Y) \f X \v_{\cal F} v_1 : X \f Y$. $v_1$ ne peut pas \^etre une variable. Donc on a de
nouveau deux cas \`a voir : 
\begin{itemize}
\item[] -- Si $v_1 = (z ~ v'_1...v'_k)$ ($k \geq 1$), alors $z = \a$, et $\G , x : (X \f Y)
\f X \v_{\cal F} v'_1 : P$. Ce qui est impossible par hypoth\`ese d'induction.
\item[] -- Si $v_1 = \l y w$, alors $\G , x : (X \f Y) \f X , y : X \v_{\cal F} w : Y$. $w$
ne peut pas \^etre ni une variable ni un $\l$-terme qui commence par $\l$. Donc $w = (z ~
w_1...w_r)$ ($r \geq 1$) et $z \not = x,y,x_i,y_j$. Donc $z = \a$, et $\G , x : (X \f Y) \f X
, y : X  \v_{\cal F} w_1 : P$. Ce qui est impossible par hypoth\`ese d'induction.
\end{itemize}
\end{itemize}
2) Si $\v_{\cal F_C} (Q \f P),P \f \q XX$, alors $\v_{\cal F_C} \q XX$ (puisque $\v_{\cal F_C}
P$). Ce qui contredit le th\'eor\`eme 5.  \hfill $\Box$ \\

{\bf Lemme 5} {\it Soit $w$ un $\l$-terme normal. Si $\a : Q \f P, x : X, f : X \f X \v_{\cal
F} w : X$, alors il existe un $n \in {\bf N}$ tel que $w = (f^n ~ x)$.} \\

{\bf Preuve} Par induction sur $w$.\\
Le $\l$-terme $w$ ne peut pas commencer par un $\l$. Si $w$ est une variable, alors $w = x$. 
Donc $w = (y ~ w_1...w_m)$ ($m \geq 1$), et on a deux possibilit\'es pour la variable $y$.
\begin{itemize}
\item[] -- Si $y = \a$, alors $\a : Q \f P, x : X, f : X \f X \v_{\cal F} w_1 : Q$, et
$Q \f P, X, X \f X \v_{\cal F} Q$. Soit $U$ une formule close d\'emontrable dans le syst\`eme
logique $\cal F$. D'apr\`es le lemme 2, on a  $\{ Q \f P, X, X \f X \} [ U / X ] \v_{\cal F} Q
[ U / X ]$, donc  $Q \f P \v_{\cal F} Q$ (puisque $U$ et $U \f U$ sont d\'emontrables). Ce qui
contredit 2) du lemme 4. 
\item[] -- Si $y = f$, alors $m = 1$ et $\a : Q \f P, x : X, f : X \f X \v_{\cal F} w_1 : X$.
Par  hypoth\`ese d'induction, il existe un $n \in {\bf N}$ tel que $w_1 = (f^n ~ x)$, donc $w
= (f^{n+1} ~ x)$. \hfill $\Box$  
\end{itemize} 

Soit $D = (Q \f P) \f N$ et, pout tout $n \in {\bf N}$, $d_n = \l \a \so{n}$. \\

{\bf Lemme 6} {\it Soit $t$ un $\l$-terme normal clos. $\v_{\cal F} t : D$ ssi il existe un 
$n \in {\bf N}$ tel que $t = d_n$.} \\

{\bf Preuve} 
$\Longleftarrow$) Facile \`a v\'erifier. \\
$\Longrightarrow$) Comme $t$ est clos, alors $t = \l \a u$ et $\a : Q \f P \v_{\cal F} u :
N$. $u$ ne peut pas \^etre une variable et si $u = (\a ~ u_1...u_m)$ ($m \geq 1$), alors $\a
: Q \f P \v_{\cal F} u_1 : Q$, donc $Q \f P \v_{\cal F} Q$. Ce qui contredit 2) du lemme 4.
Donc $u = \l x v$, et $\a : Q \f P, x : X \v_{\cal F} v :(X \f X) \f X$. $v$ ne peut pas
\^etre une variable, donc on a deux cas \`a voir. 
\begin{itemize}
\item[] -- Si $v = (y ~ v_1...v_m)$ ($m \geq 1$), alors $y = \a$, $\a : Q \f P, x : X
\v_{\cal F} v_1 : Q$, et $Q \f P, X \v_{\cal F} Q$. Soit $U$ une formule close d\'emontrable
dans le syst\`eme logique ${\cal F}$. D'apr\`es le lemme 2, on a 
$\{ Q \f P, X \} [ U / X ] \v_{\cal F} Q [ U / X ]$, donc  $Q \f P \v_{\cal F} Q$. Ce qui
contredit 2) du lemme 4. 
\item[] -- Si $v = \l f w$, alors $\a : Q \f P, x : X, f : X \f X
\v_{\cal F} w : X$. Donc, d'apr\`es le lemme 5, il existe un $n \in {\bf N}$ tel que $w = (f^n ~
x)$ et $t = d_n$. \hfill $\Box$ 
\end{itemize}

{\bf Lemme 7} Soit $S_d = \l n \l \a (\so{S} ~ (n ~ \a))$. \\
{\it 1) Pour tout $n \in {\bf N}$, $(S_d ~ d_n) \simeq_{\b} d_{n+1}$.\\
2) Pour tout $n \in {\bf N}$,  $({S_d}^n ~ d_0) \simeq_{\b} d_n$.} \\

{\bf Preuve} Facile \`a v\'erifier. \hfill $\Box$\\

{\bf Lemme 8} {\it Soit $O_d = \l n (n ~ T_P ~ \ch{d_0} ~ \ch{S_d})$ o\`u\\ $T_P = \l \a t_P$,
$\ch{d_0} = \l f (f ~ d_0)$, et $\ch{S_d} = \l x \l y (x ~ \l z (y ~ (S_d ~  z)))$.\\
$O_d$ est un op\'erateur de mise en m\'emoire pour le syst\`eme num\'erique typ\'e $<D,{\bf
d}>$.}\\

{\bf Preuve} On va d\'emontrer que :\\
1) $D \inc D$* et $\v_{\cal F} O_d : D$*$\f \neg \neg D$.\\ 
2) Pour tout $\th_n \simeq_{\b} d_n$, $(O_d ~ \th_n ~ f) \p (f ~ ({S_d}^n ~ d_0))$. \\

1) Il est facile de v\'erifier que $D \inc D$*.\\
On a :
\begin{eqnarray}
\v_{\cal F} d_0 : D &\Longrightarrow &\v_{\cal F} \ch{d_0} : \neg \neg D \nonumber 
\end{eqnarray}
et
\begin{eqnarray}
\v_{\cal F} S_d : D \f D &\Longrightarrow &y: \neg D, z: D \v_{\cal F} (y ~ (S_d ~  z)) :
\perp \nonumber \\
&\Longrightarrow &x: \neg \neg D, y: \neg D \v_{\cal F} (x ~ \l z (y ~ (S_d
~  z))) : \perp \nonumber \\
&\Longrightarrow &\v_{\cal F} \ch{S_d} : \neg \neg D \f \neg \neg D.\nonumber 
\end{eqnarray}
De plus, d'apr\`es le lemme 3, on a $\v_{\cal F} t_P : P$*, donc $\v_{\cal F} T_P : (Q \f
P)$*$=P$*$\f Q$*. \\
D'o\`u
\begin{eqnarray}
n:D^{*}\v_{\cal F} (n ~ T_P) : N^{*} &\Longrightarrow &n:D^{*}\v_{\cal F} (n ~ T_P) : \neg \neg
D, (\neg \neg D \f \neg \neg D) \f \neg \neg D \nonumber \\
&\Longrightarrow &\v_{\cal F} O_d : D^{*}\f \neg \neg D. \nonumber
\end{eqnarray}

2) Soit $\th_n \simeq_{\b} d_n$. 
\begin{itemize}
\item[] -- Si $n=0$, alors $\th_n \p d_0$.
\item[] -- Si $n \neq 0$, alors $\th_n \p \l \a \l x \l g (g ~ t_{n-1})$, $t_{n-k}
\p (g ~ t_{n-k-1})$ ($1 \leq k \leq n-1$), et $t_0 \p x$.
\end{itemize}

Si $n=0$, alors 
\begin{eqnarray}
(O_d ~ \th_n ~ f) &\sim &(\ch{d_0} ~  f) \nonumber \\ 
&\sim &(f ~ d_0). \nonumber
\end{eqnarray}

Si $n \neq 0$, alors 
\begin{eqnarray}
(O_d ~ \th_n ~ f) &\sim &(\ch{S_d} ~ t_{n-1}[\ch{S_d}/g,\ch{d_0}/x] ~ f) \nonumber \\ 
&\sim &(t_{n-1}[\ch{S_d}/g,\ch{d_0}/x] ~ \l z (f ~ (S_d ~  z))).\nonumber
\end{eqnarray}

On d\'efinit deux suites de $\l$-termes $(\t_i)_{1 \leq i \leq n}$ :
 
\begin{center}
$\t_1=\l z (f ~ (S_d ~  z))$\\
et
$\t_{k+1}=\l z (\t_k ~ (S_d ~  z))$ pour tout $(1 \leq k \leq n-1)$
\end{center} 

D\'emontrons (par r\'ecurrence sur $k$) que, pour tout $(1 \leq k \leq n)$, on a :

\begin{center}
$(O_d ~ \th_n ~ f) \sim (t_{n-k}[\ch{S_d}/g,\ch{d_0}/x] ~ \t_k)$ 
\end{center} 

\begin{itemize}
\item[] -- Pour $k=1$, le r\'esultat est vrai.
\item[] -- Supposons le r\'esultat vrai pour $k$, et d\'emontrons le pour $k+1$. 
\begin{eqnarray}
(O_d ~ \th_n ~ f) &\sim &(t_{n-k}[\ch{S_d}/g,\ch{d_0}/x] ~ \t_k) \nonumber \\ 
&\sim &(\ch{S_d} ~ t_{n-k-1}[\ch{S_d}/g,\ch{d_0}/x] ~ \t_k) \nonumber \\ 
&\sim &(t_{n-k-1}[\ch{S_d}/g,\ch{d_0}/x] ~ \l z(\t_k ~ (S_d ~ z))) \nonumber \\ 
&= &(t_{n-k-1}[\ch{S_d}/g,\ch{d_0}/x] ~ \t_{k+1}). \nonumber
\end{eqnarray}
\end{itemize}

Donc, en particulier, pour $k=n$ on a :
\begin{eqnarray}
(O_d ~ \th_n ~ f) &\sim &(t_0[\ch{S_d}/g,\ch{d_0}/x] ~ \t_k) \nonumber \\ 
&\sim &(\ch{d_0} ~ \t_n) \nonumber \\ 
&\sim &(\t_n ~ d_0).\nonumber
\end{eqnarray}

D\'emontrons (par r\'ecurrence sur $k$) que, pour tout $(1 \leq k \leq n)$, on a :

\begin{center}
$\t_k \sim \l z(f ~ ({S_d}^k ~ z))$
\end{center}

\begin{itemize}
\item[] -- Pour $k=1$, le r\'esultat est vrai.
\item[] -- Supposons le r\'esultat vrai pour $k$, et d\'emontrons le pour $k+1$. 
\begin{eqnarray}
\t_{k+1} &= &\l z(\t_k ~(S_d ~z)) \nonumber \\ 
&\sim &\l z(\l z(f~ ({S_d}^k ~ z)) ~(S_d ~ z)) \nonumber \\  
&\sim &\l z(f ~ ({S_d}^{k+1} ~ z)). \nonumber
\end{eqnarray}
\end{itemize}

Donc, en particulier, pour $k=n$ on a : $\t_n \sim \l z(f ~ ({S_d}^n ~ z))$. \\
Et
\begin{eqnarray}
(O_d ~ \th_n ~ f) &\sim &(\l z(f ~ ({S_d}^n ~ z)) ~ d_0) \nonumber \\  
&\sim &(f ~ ({S_d}^n ~ d_0)). \nonumber 
\end{eqnarray}
D'o\`u $(O_d ~ \th_n ~ f) \p (f ~ ({S_d}^n ~ d_0))$. \hfill $\Box$ \\

{\bf Lemme 9} Soit $O_B = \l n ( n~ \l f(f ~ T) ~ \l f(f ~ F))$.\\
{\it 1) $\v_{\cal F} O_B : B$*$\f \neg \neg B$. \\
2) Pour tout $\ep = T$ ou $F$ et pour tout $\th_{\ep} \simeq\sb{\b} \ep$, $(O_B ~ \th_{\ep} ~
f) \p (f ~ \ep)$.} \\ 

{\bf Preuve}  \\
1) On a :
\begin{eqnarray}
\v_{\cal F} T : B &\Longrightarrow &\v_{\cal F} \l f(f ~ T) : \neg \neg B \nonumber 
\end{eqnarray}
et
\begin{eqnarray}
\v_{\cal F} F : B &\Longrightarrow &\v_{\cal F} \l f(f ~ F) : \neg \neg B. \nonumber 
\end{eqnarray}
Donc
\begin{eqnarray}
n:B^{*}\v_{\cal F} n : \neg \neg B, \neg \neg B \f \neg \neg B &\Longrightarrow &\v_{\cal F}
O_B : B^{*}\f \neg \neg B. \nonumber 
\end{eqnarray}

2) Si $\th_{\ep} \simeq\sb{\b} \ep$, alors $\th_{\ep} \p \ep$, et donc
\begin{eqnarray}
(O_B ~ \th_{\ep} ~ f) &\sim  &(\th_{\ep} ~ \l f(f ~ T) ~ \l f(f ~ F) ~ f) \nonumber \\
&\sim &(\ep ~ \l f(f ~ T) ~ \l f(f ~ F) ~ f) \nonumber \\
&\sim  &(f ~ \ep).  \nonumber
\end{eqnarray}
D'o\`u  $(O_B ~ \th_{\ep} ~ f) \p (f ~ \ep)$. \hfill $\Box$ \\

{\bf Lemme 10} {\it Soit $t$ un $\l$-terme normal clos.  $\v_{\cal F} t : D \f B$ ssi $t
= \l \a T$ ou $t = \l \a F$.} \\

{\bf Preuve} $\Longleftarrow$) Facile \`a v\'erifier. \\
$\Longrightarrow$) Comme $t$ est clos, alors $t = \l \a u$ et $\a : D \v_{\cal F} u : B$.
$u$ ne peut pas \^etre une variable, donc on a deux cas \`a voir. 
\begin{itemize}
\item[] -- Si $u = (\a ~ u_1...u_m)$ ($m \geq 1$), alors $\a : D \v_{\cal F} u_1 : Q \f P$,
donc $\v_{\cal F} Q \f P$ (car $D$ est d\'emontrable dans le syst\`eme logique $\cal F$). Ce
qui contredit 1) du lemme 4. 
\item[] -- Si $u = \l x v$, alors $\a : D, x : X \v_{\cal F} v :X \f X$.
$v$ ne peut pas \^etre une variable, donc on a  de nouveau deux cas \`a voir. 
\begin{itemize}
\item[] -- Si $v = (z ~ v_1...v_m)$ ($m \geq 1$), alors $y = \a$, $\a : D, x : X
\v_{\cal F} v_1 : Q \ P$, et $D, X \v_{\cal F} Q \f P$. Soit $U$ une formule close
d\'emontrable dans le syst\`eme logique $\cal F$. D'apr\`es le lemme 2, on a $\{ D, X \} [ U /
X ] \v_{\cal F} \{ Q \f P \} [ U / X ]$, donc  $\v_{\cal F} Q \f P$.  Ce qui contredit 1) du
lemme 4. 
 \item[] -- Si $v = \l y w$, alors $\a : D, x : X, y : X  \v_{\cal F} w : X$.
$w$ ne peut pas commencer par un $\l$ et si $w$ est une variable, alors $w = x$ ou $w = y$,
donc $t = \l \a T$ or $t = \l \a F$. Il reste donc le cas o\`u $w = (\a ~ w_1...w_m)$ ($m
\geq 1$). Dans ce cas on a  $\a : D, x : X, y :X \v_{\cal F} w_1 : Q \f P$, et $D, X \v_{\cal
F} Q \f P$. Soit $U$ une formule close d\'emontrable dans le syst\`eme logique $\cal F$.
D'apr\`es le lemme 2, on a $\{ D, X, X \} [ U / X ] \v_{\cal F} \{ Q \f P \} [ U / X ]$,
donc  $\v_{\cal F} Q \f P$.  Ce qui contredit 1) du lemme 4.  \hfill $\Box$   
\end{itemize}
\end{itemize}

{\bf Lemme 11} Soit $t$ un $\l$-terme normal.\\
{\it 1)  Si $x :B,D \f X , y : X \v_{\cal F} t : B$, alors $t= T$ ou $t=F$. \\ 
2) Si $x :B,D \f X , y : X \v_{\cal F} t : D$, alors il existe $n \in {\bf N}$ tel que  $t =
d_n$.}\\

{\bf Preuves} M\^eme preuve que celles des lemmes 5 et 6. \hfill $\Box$ \\

Soit $E = \q X \{ ((B,D \f X),X \f X \}$. \\

Pour tous $\l$-termes $u,v$, on note $\< u,v \>$ le $\l$-terme $\l x \l y (x ~ u ~ v)$.\\

{\bf Lemme 12} {\it Soit $t$ un $\l$-terme normal clos.  $\v_{\cal F} t : E$ ssi ($t = F$)
ou il existe $n \in {\bf N}$ tel que ($t = \<b,d_n\>$  o\`u $b= T$ ou $F$).} \\

{\bf Preuve}  $\Longleftarrow$) Facile \`a v\'erifier. \\
$\Longrightarrow$) Soit $t$ un $\l$-terme normal clos tel que $\v_{\cal F} t : E$. Alors
$t = \l x u$ et $x :B,D \f X \v_{\cal F} u : X \f X$. $u$ ne peut pas \^etre une
variable, donc on a deux cas \`a voir. 
\begin{itemize}
\item[] - Si $u = (x ~ u_1...u_m)$ ($m \geq 1$), alors $x :B,D \f X \v_{\cal F} (x ~ u_1)
: D \f X$. Ce qui est impossible. 
\item[] - Si $u = \l y v$, alors $x :B,D \f X , y : X \v_{\cal F} v : X$.
$v$ ne peut pas commencer par un $\l$, donc on a de nouveau deux cas \`a voir.
\begin{itemize}
\item[] - Si $v$ est une variable, alors $v=y$ et $u = F$.
\item[] - Si $v = (z ~ v_1...v_m)$ ($m \geq 1$), alors $z = x$, $n=2$, $x :B,D \f X , y
: X \v_{\cal F} v_1 : D$, et $x :B,D \f X , y : X \v_{\cal F} v_2 : B$. Donc, d'apr\`es le
lemme 11, il existe $n \in {\bf N}$ tel que $t = \<b,d_n\>$ o\`u $b= T$ ou $b=F$.
\hfill $\Box$  \end{itemize}
\end{itemize}

{\bf Th\'eor\`eme  9} {\it Il existe un syst\`eme num\'erique typ\'e non ad\'eaquat qui
poss\`ede un op\'erateur de mise en m\'emoire.} \\

{\bf Preuve} Soit ${\cal E} = <E,{\bf e}>$ o\`u :
\begin{eqnarray} 
e_0 &= &F   \nonumber \\
e_{2n+1} &= &{\<F,d_n \>} \quad (n \geq 0) \nonumber \\
e_{2n+2} &= &{\<T,d_n \>} \quad (n \geq 0) \nonumber 
\end{eqnarray} 

\so{\so{{\bf Le test \`a z\'ero}}} \\

Soit $Z_e = \l n (n ~ \l x \l y F ~ T)$.

\begin{itemize}

\item[] \so{{\bf Typage de $Z_e$}}\\

On a :
\begin{eqnarray} 
x:B , y:D \v_{\cal F} F : B &\Longrightarrow &\v_{\cal F} \l x \l y F : B,D \f B \nonumber
\end{eqnarray}
donc
\begin{eqnarray}
n:E \v_{\cal F} n : (B,D \f B),B \f B &\Longrightarrow &n:E \v_{\cal F} (n ~ \l x \l y F ~ T) :
B \nonumber \\
&\Longrightarrow &\v_{\cal F} Z_e : E \f B. \nonumber
\end{eqnarray}

\item[] \so{{\bf Fonctionnement de $Z_e$}} \\

Si $n = 0$, alors :
\begin{eqnarray}  
(Z_e ~ e_n) &\simeq_{\b} &(F ~ \l x \l y F ~ T) \nonumber \\
&\simeq_{\b} &T. \nonumber
\end{eqnarray}
Si $n \not = 0$, alors :
\begin{eqnarray} 
(Z_e ~ e_{n+1}) &\simeq_{\b} &(e_{n+1} ~ \l x \l y F ~ T) \nonumber \\
&\simeq_{\b} &(\l x \l y F ~ b ~ d_m) \nonumber \\ 
&\simeq_{\b} &F.\nonumber
\end{eqnarray}

\end{itemize}

\so{\so{{\bf Le successeur}}} \\

Soit $S_e = \l n ((Z_e ~ n) ~ 
																											e_1 ~ 
																											((n ~ T ~ T) ~ 
																															\<F ,(S_d  ~ (n ~ F ~ d_0)) \> ~
 																														\< T , (n ~ F ~ d_0) \> ))$.

\begin{itemize}

\item[] \so{{\bf Typage de $S_e$}} \\

On a :
\begin{eqnarray} 
n: E\v_{\cal F} n : (B,D \f D), D \f D &\Longrightarrow &n: E\v_{\cal F} (n ~ F ~ d_0) : D \nonumber \\
&\Longrightarrow &n: E\v_{\cal F} \< T , (n ~ F ~ d_0) \> : E \nonumber
\end{eqnarray}
et
\begin{eqnarray} 
n: E\v_{\cal F} n : (B,D \f D), D \f D &\Longrightarrow &n: E\v_{\cal F} (S_d  ~ (n ~ F ~
d_0)) : D \nonumber \\ 
&\Longrightarrow &n: E \v_{\cal F} \<F ,(S_d  ~ (n ~ F ~ d_0)) \> : E.\nonumber
\end{eqnarray}
Donc
\begin{eqnarray}
n: E\v_{\cal F} n : (B,D \f B), B \f B &\Longrightarrow &n: E\v_{\cal F} (n ~ T ~ T) : B\nonumber \\
&\Longrightarrow &n: E\v_{\cal F} (n ~ T ~ T) : E,E \f E \nonumber \\
&\Longrightarrow &n: E\v_{\cal F} ((n ~ T ~ T) ~  \<F ,(S_d  ~ (n ~ F ~ d_0)) \> \nonumber \\
&\quad &{\< T , (n ~
F ~ d_0) \> ) : E.} \nonumber 
\end{eqnarray}
D'o\`u
\begin{eqnarray}
n: E\v_{\cal F} (Z_e ~ n) : B &\Longrightarrow &n: E\v_{\cal F} (Z_e ~ n) : E,E \f E
\nonumber \\ &\Longrightarrow &\v_{\cal F} S_e : E \f E.\nonumber
\end{eqnarray}

\item[] \so{{\bf Fonctionnement de $S_e$}} \\

On a trois cas :
\begin{eqnarray} 
(S_e ~ e_0) &\simeq_{\b} &(T ~ e_1 ~ ((e_0 ~ T ~ T) ~ \<F ,(S_d  ~ (e_0 ~ F ~
d_0)) \> ~ \< T , (e_0 ~ F ~ d_0) \> )) \nonumber \\
&\simeq_{\b} &e_1. \nonumber 
\end{eqnarray}

\begin{eqnarray} 
(S_e ~ \<F , d_n \>) &\simeq_{\b} &(F ~ e_1 ~ ((\< F , d_n \> ~ T ~ T) ~
\<F ,(S_d  ~ (e_0 ~ F ~ d_0)) \> \nonumber \\
&\quad &{\< T , (\< F , d_n \> ~ F ~ d_0) \> ))} \nonumber \\
&\simeq_{\b} &((\< F , d_n \> ~ T ~ T) ~ \<F ,(S_d  ~ (e_0 ~ F ~ d_0)) \> \nonumber \\
&\quad &{\< T , (\< F , d_n \> ~ F ~ d_0) \> )} \nonumber \\
&\simeq_{\b} &(F ~ \<F ,(S_d  ~ (e_0 ~ F ~ d_0)) \> ~ \< T , (\< F , d_n \>
~ F ~ d_0) \> )\nonumber \\
&\simeq_{\b} &{\< T , (\< F , d_n \> ~ F ~ d_0) \>} \nonumber \\
&\simeq_{\b} &{\< T , d_n \>}. \nonumber 
\end{eqnarray}
\begin{eqnarray}
(S_e ~ \< T , d_n \>) &\simeq_{\b} &(F ~ e_1 ~ ((\< T , d_n \> ~ T ~ T)
~ \<F ,(S_d  ~ (\< F , d_n \> ~ F ~ d_0)) \> \nonumber \\
&\quad &{\< T , (\< T , d_n \> ~ F ~ d_0) \> ))} \nonumber \\ 
&\simeq_{\b} &((\< T , d_n \> ~ T ~ T) ~ \<F ,(S_d  ~ (\< F , d_n \> ~ F ~ d_0)) \> \nonumber \\
&\quad &{\< T , (\< T , d_n \> ~ F ~ d_0) \> )} \nonumber \\ 
&\simeq_{\b} &(T ~ \<F ,(S_d  ~ (\< T , d_n \> ~ F ~ d_0)) \> ~ \< T , (\< T , d_n \> ~ F ~
d_0) \> ) \nonumber \\  
&\simeq_{\b} &{\<F ,(S_d  ~ (\< T , d_n \> ~ F ~ d_0)) \>} \nonumber \\
&\simeq_{\b} &{\< F , (S_d  ~ d_n) \>} \nonumber \\
&\simeq_{\b} &{\< F , d_{n+1} \>}. \nonumber   
\end{eqnarray}

\end{itemize}

\so{\so{{\bf L'op\'erateur de mise en m\'emoire}}} \\

Il est facile de v\'erifier que  $E \inc E$*.\\

Soit $O_e =\l n (n ~ \ch{S_e} ~ \ch{e_0})$ o\`u $\ch{S_e} = \l x \l y \l z ((O_B ~ x) ~ \l u
((O_d ~ y) ~ \l v (z ~ \< u,v \> )))$ et $\ch{e_0} = \l f (f ~ e_0)$.

\begin{itemize}

\item[] \so{{\bf Typage de $O_e$}} \\

On a :
\begin{eqnarray}
\v_{\cal F} e_0 : E &\Longrightarrow &\v_{\cal F} \ch{e_0} : \neg \neg E. \nonumber   
\end{eqnarray}
D'autre part, en utilisant les lemmes 8 et 9, on a :
\begin{eqnarray}
u:B, v:D \v_{\cal F} \< u,v \> : E &\Longrightarrow &u:B, z: \neg E \v_{\cal F} \l v (z ~ \<
u,v \> ) : \neg D \nonumber \\
&\Longrightarrow &y:D^{*},z: \neg E \v_{\cal F} \l u((O_d ~ y) ~ \l v (z ~
\< u,v \> )) : \neg B \nonumber \\
&\Longrightarrow &x:B^{*},y:D^{*}\v_{\cal F} \l z ((O_B ~ x) ~ \l u
((O_d ~ y) \nonumber \\
&\quad &{\l v (z ~ \< u,v \> ))) : \neg \neg E} \nonumber \\
&\Longrightarrow &\v_{\cal F}  \ch{S_e} : B^{*},D^{*}\f \neg \neg E. \nonumber   
\end{eqnarray}
D'o\`u
\begin{eqnarray}
n:E^{*}\v_{\cal F} n : (B^{*},D^{*}\f \neg \neg E) , \neg \neg E \  \neg \neg E
&\Longrightarrow &\v_{\cal F} O_e : E^{*}\f \neg \neg E. \nonumber   
\end{eqnarray}

\item[] \so{{\bf Fonctionnement de $O_e$}} \\

Soit $\th_n \simeq_{\b} e_n$, alors : 
\begin{itemize}
\item[] -- Si $n=0$, alors $\th_n \p e_0$.
\item[] -- Si $n \neq 0$, alors $\th_n \p \l x \l y (x ~ \a_n ~ \b_n)$ o\`u
$\a_n \simeq_{\b} \ep$ et $\b_n \simeq_{\b} d_m$ si $e_n = \<\ep , d_m \>$.
\end{itemize}

Si $n=0$, alors 
\begin{eqnarray}
(O_e ~ \th_n ~ f) &\sim &(e_0 ~ \ch{S_e} ~ \ch{e_0} ~ f) \nonumber \\
&\sim &(\ch{e_0} ~ f) \nonumber \\
&\sim &(f  ~ e_0). \nonumber   
\end{eqnarray}

Si $n \neq 0$, alors 
\begin{eqnarray}
(O_e ~ \th_n ~ f) &\sim &(\ch{S_e} ~ \a_n ~ \b_n ~ f) \nonumber \\
&\sim &((O_B ~ \a_n) ~ \l u ((O_d ~ \b_n) ~ \l v (f ~ \< u,v \> ))). \nonumber   
\end{eqnarray} 

D'apr\`es le Lemma 9, on a : pour tout $\l$-terme $U$, $((O_B ~ \a_n) ~ U) \sim (U ~
\ep)$.\\

Donc 
\begin{eqnarray}
(O_e ~ \th_n ~ f) &\sim &(\l u ((O_d ~ \b_n) ~ \l v (f ~ \< u,v \> ))  ~ \ep) \nonumber \\
&\sim &((O_d ~ \b_n) ~ \l v (f ~ \<\ep,v \> )). \nonumber   
\end{eqnarray} 
D'apr\`es le Lemma 8, on a : pour tout $\l$-terme $V$, $((O_d ~ \b_n) ~ V) \sim (V ~
({S_d}^m ~ d_0))$.\\

Donc  
\begin{eqnarray}
(O_e ~ \th_n ~ f) &\sim &(\l v (f ~ \< \ep,v \> ) ~ ({S_d}^m ~ d_0)) \nonumber \\
&\sim &(f ~ \< \ep,({S_d}^m ~ d_0) \> ). \nonumber   
\end{eqnarray} 
D'o\`u $(O_e ~ \th_n ~ f) \p (f ~ \< \ep,({S_d}^m ~ d_0) \> )$
\end{itemize}

\so{\so{{\bf L'in\'existance d'un pr\'ed\'ecesseur}}} \\

Supposons qu'il existe un $\l$-terme normal clos $P_e$ pour le pr\'ed\'ecesseur. \\
Soit $P' = \l n (P_e ~ \<F , n \> ~ T ~ F)$. \\

On a
\begin{eqnarray}
n:D \v_{\cal F} \<F , n\> : E &\Longrightarrow &n:D \v_{\cal F} (P_e ~ \<F , n \>) : E
\nonumber \\
&\Longrightarrow &n:D \v_{\cal F} (P_e ~ \<F , n \>) : (B,D \f B),B \f B \nonumber \\
&\Longrightarrow &n:D \v_{\cal F} (P_e ~ \<F , n \> ~ T ~ F) : B \nonumber \\
&\Longrightarrow &\v_{\cal F} P' : D \f B. \nonumber 
\end{eqnarray} 
Donc, d'apr\`es le lemme 10, on obtient $P' = \l \a T$ ou $P' = \l \a F$.\\
Mais on a : 
\begin{eqnarray}
(P' ~ d_0) &\simeq_{\b} &(P_e ~ e_1 ~ T ~ F) \nonumber \\
&\simeq_{\b} &(e_0 ~ T ~ F) \nonumber \\
&\simeq_{\b} &F \nonumber  
\end{eqnarray} 
et
\begin{eqnarray}
(P' ~ d_1) &\simeq_{\b} &(P_e ~ e_3 ~ T ~ F) \nonumber \\
&\simeq_{\b} &(e_2 ~ T ~ F) \nonumber \\
&\simeq_{\b} &(T ~ T ~ d_0) \nonumber \\
&\simeq_{\b} &T.\nonumber 
\end{eqnarray} 
D'o\`u une contradiction. \hfill $\Box$\\

{\bf Remarque} : Il est facile de v\'erifier que le $\l$-terme \\
$P_e = \l n ((Z_e ~ n) ~ 
																				e_0 ~ 
																									((n ~ T ~ T) ~
																															\< F , (n ~ F ~ d_0) \>  ~
																																		(\so{Z} ~ (n ~ F ~ d_0 ~ T)) ~
																															            \<T ,(\l \a (\so{P}  ~ (n ~ F ~ d_0 ~ \a))) \> ~
																																            F)
																																																																																									)
																														                      																																					)$
est un pr\'ed\'ecesseur (non typable dans le syst\`eme ${\cal F}$ de type $E \f E$) pour le
syst\`eme num\'erique {\bf e}.\\

\section{Conclusion}

Suite \`a cette \'etude, deux questions se posent :
\begin{itemize}
\item Est-il vrai que chaque syst\`eme num\'erique typ\'e ad\'equat poss\`ede un op\'erateur
de mise en m\'emoire? En effet l'op\'erateur de mise en m\'emoire qu'on a construit pour un
syst\`eme num\'erique ad\'equat quelconque (voir la preuve du th\'eor\`eme 6) utilise un
op\'erateur de point fixe et donc il est non typable dans le syst\`eme $\cal F$.  
\item Quelles sont les fonctions qu'on peut repr\'esenter dans un syst\`eme num\'erique
typ\'e ad\'equat? 
\end{itemize}


\begin{thebibliography}{99}

\bibitem{K} H. Barendregt.
{\em  The lambda calculus, its syntax and semantics.} \\
{\bf North Holland, 1984}

\bibitem{K} J.-Y. Girard, Y. Lafont, P. Taylor.
{\em Proofs and types.} \\
{\bf 	Cambridge University Press, 1986.}

\bibitem{K} J.-L. Krivine.
{\em  Lambda calcul, types et mod\`eles.} \\
{\bf Masson, 1990}

\bibitem{K} J.-L. Krivine.
{\em Op\'erateurs de mise en m\'emoire et traduction de G\H{o}del}\\
{\bf Archive for Mathematical Logic 30 (1990), pp. 241-267.}

\bibitem{Z} K. Nour.
{\em Op\'erateurs de mise en m\'emoire en lambda-calcul pure et typ\'e}\\
{\bf Th\`ese de Doctorat, Universit\'e de Chamb\'ery, 1993.}

\bibitem{K} K. Nour.
{\em An example of a non adequate numeral system.} \\
{\bf CRAS. Paris, 323, S\'erie I (1996), pp. 439-442.}

\bibitem{K} K. Nour.
{\em A conjecture on numeral system.} \\
{\bf Notre Dame of Formal Logic, vol. 38 (1997), pp. 270-275.}

\end{thebibliography}
\end{document}